# AI-Driven Adaptive Air Transit Network with Modular Aerial Pods

Amir Shafiee[a,1], Alireza Yazdiani[b,1], Hanieh Rastegar[a,2], Rui Li[b,2], Rayan Karim[b], Aolei Cao[b], Ziyang Li[b], Xieqing Yu[b], Charlle Sy[b], Zhaoyao Bao[b], Xi Cheng[b,*], H. Oliver Gao[b]

[a]*Department of Civil, Material and Environmental Engineering, The University of Illinois Chicago, Chicago, IL 60607, USA*
[b] *Systems Engineering, Cornell University, Ithaca, NY 14850, USA*

**Abstract**

This paper presents an adaptive air transit network leveraging modular aerial pods and artificial intelligence (AI) to address urban mobility challenges. Passenger demand, forecasted from AI models, serves as input parameters for a Mixed-Integer Nonlinear Programming (MINLP) optimization model that dynamically adjusts pod dispatch schedules and train lengths in response to demand variations. The results reveal a complex interplay of factors, including demand levels, headway bounds, train configurations, and fleet sizes, which collectively influence network performance and service quality. The proposed system demonstrates the importance of dynamic adjustments, where modularity mitigates capacity bottlenecks and improves operational efficiency. Additionally, the framework enhances energy efficiency and optimizes resource utilization through flexible and adaptive scheduling. This framework provides a foundation for a responsive and sustainable urban air mobility solution, supporting the shift from static planning to agile, data-driven operations.

## 1. Introduction

As urban populations continue to grow and cities expand, the limitations of traditional ground and underground transportation systems are becoming increasingly evident. Traffic congestion, overcrowded subways, and the environmental impact of fossil fuel-based vehicles highlight the urgent need for innovative and sustainable solutions. Urban air mobility (UAM) has emerged as a promising alternative, offering the potential to alleviate these challenges by introducing a third dimension to urban transportation. By utilizing the airspace above cities, UAM can significantly reduce travel times, improve connectivity, and enhance the overall quality of life in densely populated areas.

The development of advanced UAM systems represents a pivotal evolution in addressing the complex challenges of urban transportation. Existing research in ground transit network design and scheduling highlights the importance of balancing operational efficiency and passenger service quality through optimized routing, vehicle scheduling, and frequency adjustments. Studies (e.g., [1], [2], [3]) highlight the importance of robust scheduling frameworks in improving transit reliability and effectively managing disruptions. Similarly, research on modular and autonomous vehicle technologies (e.g., [4], [5], [6]), emphasizes the benefits of scalability and responsiveness in contemporary transit systems. Recent studies on UAM [7], [8], highlight the potential of aerial transit to address urban mobility challenges by integrating advanced optimization techniques like Mixed-Integer Programming (MIP) with demand forecasting models. These works underscore the critical need for system-wide planning, equitable access, and operational cost management to ensure UAM's long-term feasibility and acceptance.

In response to these challenges, we propose an innovative air transit system that combines modular aerial pod technologies with AI-driven optimization. Each modular pod can operate either independently or as part of a train, offering scalable capacity and flexible scheduling possibilities. By analyzing historical and real-time passenger data, the AI framework generates accurate forecasts of passenger demand, which can be incorporated as parameters into the design and operation of the air transit network, enabling dynamic adjustments to dispatch schedules and pod train sizes based on real-time demand to minimize passenger waiting times while maintaining operational efficiency.

---

\* Corresponding author. Email: xc557@cornell.edu.
[1]Co-first authors of the article. [2]Co-second authors of the article.

To model and optimize the network, we developed a mathematical framework inspired by principles of scheduling and allocation. Our framework employs Mixed-Integer Non-Linear Programming (MINLP) to optimize pod dispatching, train lengths, and scheduling. This integrated approach ensures that the system remains responsive to fluctuating urban transportation needs and optimizes resource utilization effectively.

The remaining sections of this paper explore the system architecture, air transit network design, and results in detail, showcasing the feasibility and effectiveness of the proposed air transit network.

## 2. AI-Driven Adaptive Air Transit Network Model

   *1. Framework*

The AI-driven adaptive air transit network, as illustrated in Figure 1, forecasts passenger demand with machine learning models using historical and real-time data. These forecasts provide input parameters for a MINLP-based optimization model, which dynamically determines key operational parameters, including dispatch schedules and pod train sizes. Modular aerial pods, capable of operating independently or as part of a train, adapt to varying demand levels by optimizing fleet operations. During peak hours, pods form trains to improve energy efficiency from a mechanism similar to truck platooning, while in low-demand periods, they operate individually to reduce costs and increase flexibility. By continuously integrating data and optimizing operations, the system ensures efficient resource utilization and reduced passenger waiting times, responding effectively to urban mobility dynamics.

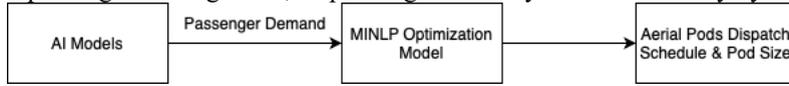

Figure 1. Overall Framework

   *2. Mathematical MINLP Model*

With the predicted passenger demand from AI models, we develop an innovative adaptive air transit system designed to operate within a grid network embedded in a densely populated urban environment, as illustrated in Figure 2, inspired from [4]. The system utilizes autonomous aerial pods that travel along predefined routes, operating across the urban buildings. There are $N \times N$ stops in the grid network in Figure 2. Layer 5 represents the road network, while Layers 1 - 4 represent the air transit network, where each layer is dedicated to one of the four cardinal directions—Westbound (WB), Eastbound (EB), Northbound (NB), and Southbound (SB), respectively. Located at the periphery of the grid (A, B, C, D), serving as origin and destination points for the aerial pods and facilitating battery swapping and maintenance operations. The system adapts to real-time passenger demand by dynamically adjusting aerial pod train lengths and departure schedules to improve service quality. When pods operate as a train, they achieve energy savings due to aerodynamic advantages, as demonstrated in [9]. While energy savings are not explicitly optimized in this model, this operational characteristic is noted in the assumptions for potential future extensions.

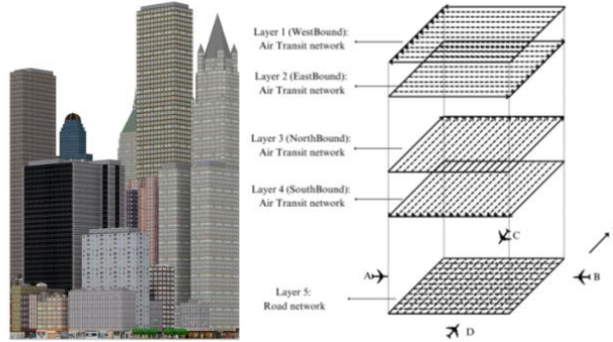

Figure 2. Conceptual diagram of the four-layer air transit network with aerial pods operating in dedicated directions.

The following assumptions are used in the air transit model:
   i. The air transit network operates on four distinct layers, each dedicated to a single direction (NB, SB, EB, WB) to prevent collisions.
   ii. Stops are integrated with urban buildings and are located at their respective layers.

iii. The network forms a grid with N lines and N stops in each direction.
iv. Pods operate in a single direction on each layer, moving between two opposite depots.
v. Pods travel at a constant cruise speed v, which includes acceleration and deceleration phases.
vi. Pods are autonomous, electric, and utilize battery swapping at depots, with a fixed swapping time.
vii. Pods can form convoys (pod trains) with a maximum allowable length $L_{max}$
viii. The length of a pod train remains constant throughout its trip and is determined based on passenger demand.
ix. Passenger arrivals at stops follow a Poisson process with rate $\lambda_{ti}$, which may vary over time.
x. Passengers are allowed at most one transfer to complete their trips.
xi. All passengers must be seated during transit; standing is not allowed.
xii. Pods dwell at each stop for a fixed time $\tau_{dwell}$
xiii. Energy savings are achieved through pod train formation, modelled as a function of pod train length.
xiv. The scheduling of pod trains aims to balance passenger service levels with operational efficiency.

Table 1: Model Parameters and Variables

| Notation | Description |
|---|---|
| **Sets and Indices** | |
| $L$ | Set of lines (indexed by $l$) |
| $D$ | Set of directions (indexed by $d$) |
| $T$ | Set of time periods (indexed by $t$) |
| $Nodes$ | Set of stations/nodes (indexed by $(i,j)$) |
| $K$ | Set of dispatch indices (indexed by $k$) |
| $OD_t$ | Set of origin-destination (OD) pairs $i,j,l,t$ representing travel demand in time period $t$ |
| **Parameters** | |
| $C$ | Capacity of a single pod (number of passengers) |
| $F$ | Total available fleet size (number of pods) |
| $h_{min}$ | Minimum allowable headway between consecutive dispatches |
| $h_{max}$ | Maximum allowable headway between consecutive dispatches |
| $T_{start}$ | Earliest start time for the first dispatch in each period |
| $T_{end}$ | Latest allowable dispatch time |
| $OD_{i,j,l,t}$ | Passenger demand for the OD pair $(i,j)$ on line $l$ during time period $t$ |
| $BigM$ | A sufficiently large constant used as a penalty for unmet passenger departures after the last dispatch |
| $Tr_l$ | Round-trip time for line $l$ (including any dwell or turnaround times) |
| **Decision Variables** | |
| $x_{l,d,t,k}$ | Departure time of the $k_{th}$ dispatch on line $l$, direction $d$, during time period $t$ |
| $L^t_{d,l,k}$ | Number of pods (train length) in the $k_{th}$ dispatch on line $l$, direction $d$, during time period $t$ |
| **State Variables** | |
| $h_{l,d,t,k}$ | Headway between dispatch $k$ and dispatch $k-1$ on line $l$, direction $d$, in time period $t$ |
| $w_{i,j,l,d,t,k}$ | Number of passengers still waiting for the $k_{th}$ dispatch with OD pair $(i,j,l,t)$ in direction $d$ |
| $s_{i,j,l,d,t,k}$ | Number of passengers seated on the $k_{th}$ dispatch with OD pair $(i,j,l,t)$ in direction $d$ |
| $OB_{i,l,d,t,k}$ | Number of onboard passengers at node on line $l$, direction $d$, for the $k_{th}$ dispatch in time period $t$ |

The notations are listed in Table 1.

$$\min Z = \sum_{(i,j,l,t)\in OD_t} \sum_{d\in D} \sum_{k\in K} w_{i,j,l,d,t,k} \cdot \begin{cases} h_{l,d,t,k+1}, & \text{if } k<|K|, \\ BigM, & \text{if } k=|K| \end{cases}$$
$$+ \sum_{(i,j,l,t)\in OD_t} \sum_{d\in D} \sum_{k\in K} \left(w_{i,j,l,d,t,k} + s_{i,j,l,d,t,k}\right) \cdot \begin{cases} \dfrac{h_{l,d,t,k}}{2}, & \text{if } k<|K|, \\ BigM, & \text{if } k=|K| \end{cases} \quad (1)$$

Subject to:

$$h_{l,d,t,k} = x_{l,d,t,k} - x_{l,d,t,k-1}, \quad \forall l \in L, d \in D, t \in T, k = 2,\cdots,|K| \tag{2}$$

$$x_{l,d,t,1} \geq T_{\text{start}}, \quad \forall l \in L, d \in D, t \in T \tag{3}$$

$$x_{l,d,t,k} \geq x_{l,d,t,k-1} + h_{\min}, \quad \forall l \in L, d \in D, t \in T, k = 2,\cdots,|K| \tag{4}$$

$$x_{l,d,t,k} \leq T_{\text{end}}, \quad \forall l \in L, d \in D, t \in T, k \in K \tag{5}$$

$$h_{l,d,t,k} \leq h_{\max}, \quad \forall l \in L, d \in D, t \in T, k \in K \tag{6}$$

$$s_{i,j,l,d_1,t,k} = 0, \quad \forall i,j \text{ such that } j \leq i, l \in L, d \in D, t \in T, k \in K \tag{7}$$

$$w_{i,j,l,d_1,t,k} = 0, \quad \forall i,j \text{ such that } j \leq i, l \in L, d \in D, t \in T, k \in K \tag{8}$$

$$s_{i,j,l,d_2,t,k} = 0, \quad \forall i,j \text{ such that } j \leq i, l \in L, d \in D, t \in T, k \in K \tag{9}$$

$$w_{i,j,l,d_2,t,k} = 0, \quad \forall i,j \text{ such that } j \leq i, l \in L, d \in D, t \in T, k \in K \tag{10}$$

$$\sum_{d \in D} \left( s_{i,j,l,d,t,1} + w_{i,j,l,d,t,1} \right) = \text{OD}_{i,j,l,t}, \quad \forall (i,j,l,t) \in \text{OD}_t \tag{11}$$

$$\sum_{d \in D} \left( s_{i,j,l,d,t,k} + w_{i,j,l,d,t,k} \right) = \text{OD}_{i,j,l,t} - \sum_{d \in D} \sum_{m \in K : m < k} s_{i,j,l,d,t,m}, \quad \forall (i,j,l,t) \in \text{OD}_t, k = 2,\cdots,|K| \tag{12}$$

$$\text{OB}_{i,l,d_1,t,k} = \begin{cases} \sum_{j>i} s_{i,j,l,d_1,t,k}, & \text{if } i = \min(\text{nodes}), \\ \text{OB}_{i-1,l,d_1,t,k} + \sum_{j>i} s_{i,j,l,d_1,t,k} - \sum_{j<i} s_{j,i,l,d_1,t,k}, & \text{otherwise} \end{cases}, \forall l \in L, t \in T, k \in K \tag{13}$$

$$\text{OB}_{i,l,d_2,t,k} = \begin{cases} \sum_{j<i} s_{i,j,l,d_2,t,k}, & \text{if } i = \max(\text{nodes}), \\ \text{OB}_{i+1,l,d_2,t,k} + \sum_{j<i} s_{i,j,l,d_2,t,k} - \sum_{j>i} s_{j,i,l,d_2,t,k}, & \text{otherwise} \end{cases}, \forall l \in L, t \in T, k \in K \tag{14}$$

$$\text{OB}_{i,l,d,t,k} \leq C \cdot L^t_{l,d,k}, \quad \forall i \in \text{nodes}, l \in L, d \in D, t \in T, k \in K \tag{15}$$

$$\sum_{l \in L} \sum_{d \in D} \sum_{k \in K} \left( \frac{Tr_l}{h_{l,d,t,k}} \cdot L^t_{l,d,k} \right) \leq F, \quad \forall t \in T \tag{16}$$

$$h_{l,d,t,k}, x_{l,d,t,k}, s_{i,j,l,d,t,k}, w_{i,j,l,d,t,k}, \text{OB}_{i,l,d,t,k}, L^t_{l,d,k} \geq 0, \quad \forall i,j \in \text{nodes}, l \in L, d \in D, t \in T, k \in K \tag{17}$$

The objective function (1) aims to minimize the total waiting time for passengers across all dispatches. This is achieved through two components. The first term minimizes the waiting time ($w_{i,j,l,d,t,k} \cdot h_{l,d,t,k+1}$) for passengers who missed the previous dispatch and must wait until the next dispatch arrives. To ensure no passengers remain waiting after the last dispatch, a penalty cost ($BigM$) is applied to any remaining waiting passengers. The second term accounts for the waiting time of all passengers arriving at the station for the current dispatch. This includes both passengers who will board the pod and those who will miss the upcoming dispatch. To approximate this waiting time, the average wait time is calculated as half of the headway (($w_{i,j,l,d,t,k} + s_{i,j,l,d,t,k}) \times \frac{h_{l,d,t,k}}{2}$). This objective ensures efficient scheduling by minimizing delays for waiting passengers while accounting for the dynamics of passenger boarding and waiting, ensuring that all passengers are served with minimal delay.

The first constraint (2) defines the headway as the difference in dispatch times between consecutive pods for each line, direction, time period, and pod index, ensuring consistent calculation of the time intervals. Constraint (3) ensures that the first dispatch for each line, direction, and time period starts no earlier than the predetermined start time ($T_{start}$). The minimum headway constraint (4) ensures a safe and feasible separation between dispatches, and the departure time limit constraint (5) keeps all dispatches within the operational timeframe ($T_{end}$). Constraint (6) ensures that the headway time does not exceed a predefined upper limit, thereby maintaining a minimum frequency of service. Directional feasibility constraints (7)-(10) govern passenger boarding to ensure passengers only board pods traveling in the correct direction. The next two constraints (11) and (12) guarantee that the total number of seated and waiting passengers matches the demand for each origin-destination pair. Constraints (13) and (14) define the recursive updates for onboard passengers in both forward ($d_1$) and reverse ($d_2$) directions. For the first or last node (depending on the direction), the onboard passengers are calculated as the sum of passengers boarding at that node for valid destinations.

For all other nodes, the onboard count is updated recursively by adding passengers boarding at the current node and subtracting those alighting, ensuring accurate tracking of onboard passenger numbers throughout the route. The capacity constraint (15) ensures that the number of onboard passengers does not exceed the pod's capacity, scaled by the total number of pods in the convoy. Constraint (16) ensures that the total number of train units deployed across all lines, directions, and time periods does not exceed the available fleet size $F$. Lastly, constraint (17) ensures that all decision variables, such as headway, dispatch times, passenger boarding, waiting, onboard counts, and train lengths, are non-negative, reflecting their real-world interpretation and feasibility.

## 3. Linearization of the Nonlinear Fleet Size Constraint

In our original model formulation, the fleet size constraint involves a nonlinear term. Specifically, consider the following nonlinear constraint:

$$\sum_{l \in L} \sum_{d \in D} \sum_{k \in K} \left( \frac{Tr_l}{h_{l,d,t,k}} \cdot L^t_{l,d,k} \right) \leq F, \quad \forall t \in T \tag{18}$$

Constraint (18) is linearized to ensure that our model can be solved using standard Mixed-Integer Linear Programming (MILP) solvers. The key idea is to discretize both $h_{l,d,t,k}$ and $L^t_{d,l,k}$, and introduce a set of binary variables to select the corresponding pair of values. This approach transforms the original nonlinear term into a linear combination of binary variables. We assume that the headway $h_{l,d,t,k}$ can only take values from a finite set $H = \{h_{min}, h_2, …, h_{max}\}$ and the number of pods per train $L^t_{d,l,k}$ is similarly restricted to $L_{train} = \{l^1, l^2, …, l^{Lmax}\}$. By discretizing these variables, we approximate the continuous solution space with a finite, though potentially fine-grained, set of options. Then, we introduce a binary variable $\delta_{l,d,t,k,hv,lv}$ for each line $l$, direction $d$, time period $t$, dispatch index $k$, headway value $hv \in H$, and train length value $lv \in L$. The variable $\delta_{l,d,t,k,hv,lv}$ equals 1 if and only if $h_{l,d,t,k} = hv$ and $L^t_{d,l,k} = lv$; otherwise, it is 0. We enforce the selection of exactly one pair $(hv, lv)$ for each $(l, d, t, k)$ via:

$$\sum_{hv \in H} \sum_{lv \in L_{train}} \delta_{l,d,t,k,hv,lv} = 1, \quad \forall l \in L, d \in D, t \in T, k \in K \tag{19}$$

We ensure that the chosen pair $(hv, lv)$ corresponds to the continuous variables $h_{l,d,t,k}$ and $L^t_{d,l,k}$ as follows:

$$h_{l,d,t,k} = \sum_{hv \in H} \sum_{lv \in L_{train}} hv \cdot \delta_{l,d,t,k,hv,lv} \tag{20}$$

$$L^t_{l,d,k} = \sum_{hv \in H} \sum_{lv \in L_{train}} lv \cdot \delta_{l,d,t,k,hv,lv} \tag{21}$$

The nonlinear term $\frac{Tr_l}{h_{l,d,t,k}} \times L^t_{d,l,k}$ can now be expressed as:

$$z_{l,d,t,k} = \sum_{hv \in H} \sum_{lv \in L_{train}} \left( \frac{Tr_l}{hv} \cdot lv \right) \delta_{l,d,t,k,hv,lv} \tag{22}$$

Since $Tr_l$, $hv$, and $lv$ are all constants (parameters or discretized sets), and $\delta_{l,d,t,k,hv,lv}$ is a binary variable, the expression for $z_{l,d,t,k}$ is linear. Therefore, we replace the original nonlinear constraint (18) with the following linear constraint:

$$\sum_{l \in L} \sum_{d \in D} \sum_{k \in K} z_{l,d,t,k} \leq F, \quad \forall t \in T \tag{23}$$

Here, the left-hand side is now a summation of linear terms $z_{l,d,t,k}$, each defined by a linear combination of binary variables $\delta_{l,d,t,k,hv,lv}$. This approach ensures the entire optimization model remains within a linear integer programming framework, facilitating the use of powerful MILP solvers and guaranteeing global optimality within the discretized solution space.

## 4. Numerical Experiments

In this section, we present the results obtained from the proposed air transit optimization model. We first outline the experimental setup and data used, then demonstrate the performance of the model under various operational scenarios. We subsequently discuss the impact of key parameters such as headway constraints, train length, and passenger demand on key performance metrics. Lastly, we provide a sensitivity analysis and compare the solution quality and computational performance across different instances.

　　1. *Experimental Setup*

We implemented the model using Gurobi 9.0 [10] on a MacBook Pro equipped with an Apple M2 Pro chip and a 10-core CPU. The network was designed as a grid with N×N stops, where N=8. Passenger arrivals were modeled using Poisson distributions with rates varying by time of day and location. Prior to solving the optimization model, a preprocessing step was conducted to reconcile node-level demand data with the line-level structure of the model. Initially, origin-destination (OD) demand was specified for each pair of nodes in the network. However, as the formulation operates at the line level—where each line represents a unidirectional path through multiple stops—the node-level OD demands were aggregated along each line. Specifically, for each line, OD pairs falling along that line were summed and assigned as aggregated demand to the corresponding line-level OD pair. This process streamlined the input data and ensured the optimization problem remained computationally tractable without distorting the underlying demand patterns.

The model also allows passengers to make at most one transfer between lines. In practice, this means that a trip spanning multiple nodes in different lines is decomposed into one or two segments: a direct segment along a single line if the OD pair lies on the same corridor, or two segments if a transfer is required. By aggregating demands at the line level, the model retains the ability to represent such transfers while avoiding the computational complexity of node-by-node flow tracking. Key parameters, including pod capacity, headway limits, and total fleet size, were incorporated into the model to define its operational boundaries and guide the optimization process.

Table 2. Key Parameter Values Used in the Experiments.

| Parameter | Symbol | Value(s) |
|---|---|---|
| Number of lines | $|L|$ | 8 |
| Directions | $|D|$ | 2 |
| Time periods | $|T|$ | 4 |
| Dispatch indices | $|K|$ | 5 |
| Capacity per pod | $C$ | 25 |
| Min/Max train length | $L_{min}, L_{max}$ | 1 / 5 pods |
| Min/Max headway | $h_{min}, h_{max}$ | 5 / 15 (minutes) |
| Fleet size | $FL$ | 1000 pods |
| Rounded trip time range | $Tr_l$ | 20-40 (minutes) |

2. *Baseline Scenario*

We begin by examining a baseline scenario characterized by average demand levels and parameters from Table 1. The objective function seeks to minimize passenger waiting time, ensuring a balance between service frequency and capacity utilization. Key results for the baseline scenario are presented in Table 3, including total waiting time, average headway, total number of pods dispatched, and overall load factor.

3. *Varying Demand Scenarios*

To assess the robustness of the system under fluctuating demand, we vary the average peak demand rates ($\lambda_i^t$) by $\pm 20\%$ $and$ $\pm 40\%$ for peak and off-peak hours and different lines relative to the baseline. shows how different demand intensities influence the average headway, and the total number of pods dispatched. Table 4 summarizes the results, for instance, increasing peak hour arrivals leads to longer headways and more pods deployed, while reducing demand allows for shorter headways and fewer pods, reflecting the model's adaptive resource allocation.

Table 3. Baseline Scenario Results.

| Metric | Value | Unit | Interpretation |
|---|---|---|---|
| Average Headway | 8.30 | Time Units | Service frequency |
| Total Pods Dispatched | 723 | Pods | Resource utilization |
| Average Load Factor | 24% | Percentage | Efficiency in capacity usage |

Table 4. Effects of Varying Demand

| Avg. Peak Hour demand | Avg. Off Peak Hour demand | Avg. Headway | Total Pods Dispatched |
|---|---|---|---|
| (7, 28) | (5, 20) | 8.76 | 773 |
| (6, 24) | (5, 20) | 8.16 | 736 |
| (5, 20) | (5, 20) | 8.3 | 723 |
| (5, 20) | (4, 16) | 5.85 | 690 |
| (5, 20) | (3, 12) | 6.75 | 641 |

### 4. Impact of Headway Constraints

We next investigate scenarios where the minimum and maximum allowable headways are adjusted. Table 5 demonstrates how tightening the headway bounds affects overall performance. Increasing the minimum headway reduces the permissible service frequency and might necessitate deploying more pods to mitigate increased waiting times. Conversely, relaxing the maximum headway allows for less frequent service but risks higher waiting times if not balanced correctly.

Table 5. Effects of Changing Headway Constraints

| Scenario | $hmin$ | $hmax$ | Avg. Headway | Total Pods Dispatched |
|---|---|---|---|---|
| Baseline | 5 | 15 | 8.3 | 723 |
| Tighter | 8 | 12 | 8.04 | 780 |
| Looser | 3 | 20 | 13.46 | 715 |

### 5. Train Length Variations

We investigate how allowing larger trains ($L_{max}$= 5, 10, 15) influences energy use proxies (as captured by z) and passenger waiting times. Adjusting the maximum train length ($L_{max}$) changes the capacity configuration. Table 6 shows that increasing $Lmax$ generally leads to a higher number of dispatched pods due to flexible capacity allocation. However, it can also influence headway and waiting times, as larger trains potentially reduce the frequency of dispatches needed.

### 6. Fleet Size Sensitivity

Finally, we examine the impact of varying the total fleet size ($FL$) from 500 to 1300 pods.

| $L_{max}$ | Avg. Headway | Total Pods Dispatched |
|---|---|---|
| Baseline | 8.3 | 723 |
| Tighter | 9.08 | 770 |
| Looser | 10.49 | 788 |

Table 7 shows that while increasing the fleet size initially improves headways and can reduce waiting times, beyond a certain threshold these benefits plateau. This indicates that simply adding more pods does not indefinitely enhance system performance.

## 5. Conclusions

We present an adaptive air transit network utilizing modular aerial pods and AI to address urban mobility challenges. Passenger demand, forecasted from AI models, serves as input parameters for a Mixed-Integer Nonlinear Programming (MINLP) optimization model. The system dynamically adjusts pod dispatch schedules and train lengths in response to demand variations. The results highlight a complex interplay of factors—demand levels, headway bounds, train configurations, and fleet sizes—that shape the network's performance and service quality. The model underscores the importance of dynamic adjustments, as heightened demand during peak hours necessitates frequent dispatches or longer trains, while modularity alleviates capacity bottlenecks. Even minor changes to headway or train lengths can significantly impact passenger waiting times, fleet utilization, and operational efficiency, emphasizing the need for flexibility and adaptability across infrastructure and scheduling.

Table 6. Sensitivity to Maximum Train Length

| $L_{max}$ | Avg. Headway | Total Pods Dispatched |
|---|---|---|
| Baseline | 8.3 | 723 |
| Tighter | 9.08 | 770 |
| Looser | 10.49 | 788 |

Table 7. Sensitivity to Maximum Fleet Size at each Time period

| Max Fleet Size | Avg. Headway | Total Pods Dispatched | Seat Utilization Index |
|---|---|---|---|
| 500 | 12.59 | 684 | 28 |
| 600 | 11.83 | 694 | 26 |
| 700 | 11.49 | 704 | 25 |
| 800 | 10.9 | 716 | 24 |
| 900 | 9.38 | 716 | 24 |
| 1000 | 8.3 | 723 | 24 |
| 1100 | 6.26 | 733 | 23 |
| 1200 | 5.14 | 741 | 23 |
| 1300 | 5.01 | 755 | 23 |

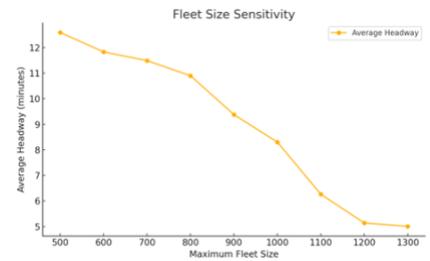

Figure 3 Impact of increasing fleet size on Avg. Headway

However, increasing fleet size alone results in diminishing returns beyond certain thresholds, underscoring the importance of strategic resource allocation. Complementary strategies, such as advanced scheduling algorithms, demand-responsive policies, and improved line synchronization, are essential to address performance gaps. Balancing operational efficiency and service quality requires the integration of user-centric metrics into decision-making processes. These findings establish a foundation for advancing adaptive air transit systems, promoting a shift from static planning to agile, data-driven, and sustainable mobility solutions.

As a proof-of-concept study, the model and demand scenarios are hypothetical, demonstrating the feasibility of the proposed concept. Future studies should validate the framework's applicability and performance in realistic operational settings, incorporating energy optimization and integration with multimodal urban mobility systems.

**References**

[1] V. Guihaire and J.-K. Hao, "Transit network design and scheduling: A global review," *Transp. Res. Part Policy Pract.*, vol. 42, no. 10, pp. 1251–1273, Dec. 2008, doi: 10.1016/j.tra.2008.03.011.
[2] H. Badia, M. Estrada, and F. Robusté, "Competitive transit network design in cities with radial street patterns," *Transp. Res. Part B Methodol.*, vol. 59, pp. 161–181, Jan. 2014, doi: 10.1016/j.trb.2013.11.006.
[3] S. Bafandkar, Y. Shafahi, A. Eslami, and A. Yazdiani, "Digitalizing railway operations: An optimization-based train rescheduling model for urban and interurban disrupted networks," *Digit. Eng.*, vol. 5, p. 100033, Mar. 2025, doi: 10.1016/j.dte.2024.100033.
[4] P. W. Chen and Y. M. Nie, "Analysis of an idealized system of demand adaptive paired-line hybrid transit," *Transp. Res. Part B Methodol.*, vol. 102, pp. 38–54, Aug. 2017, doi: 10.1016/j.trb.2017.05.004.
[5] X. Cheng, Y. (Marco) Nie, and J. Lin, "An Autonomous Modular Public Transit service," *Transp. Res. Part C Emerg. Technol.*, vol. 168, p. 104746, Nov. 2024, doi: 10.1016/j.trc.2024.104746.
[6] A. Shafiee, H. R. Moghaddam, and J. Lin, "Using Autonomous Modular Vehicle Technology as an Alternative for Last-Mile Delivery," in *2024 Forum for Innovative Sustainable Transportation Systems (FISTS)*, Feb. 2024, pp. 1–6. doi: 10.1109/FISTS60717.2024.10485532.
[7] H. Shon and J. Lee, "An optimization framework for urban air mobility (UAM) planning and operations," *J. Air Transp. Manag.*, vol. 124, p. 102720, Apr. 2025, doi: 10.1016/j.jairtraman.2024.102720.
[8] A. Straubinger, E. T. Verhoef, and H. L. F. de Groot, "Will urban air mobility fly? The efficiency and distributional impacts of UAM in different urban spatial structures," *Transp. Res. Part C Emerg. Technol.*, vol. 127, p. 103124, Jun. 2021, doi: 10.1016/j.trc.2021.103124.
[9] H. Weimerskirch, J. Martin, Y. Clerquin, P. Alexandre, and S. Jiraskova, "Energy saving in flight formation," *Nature*, vol. 413, no. 6857, pp. 697–698, Oct. 2001, doi: 10.1038/35099670.
[10] Gurobi Optimization, *Gurobi 9.0*. [Online]. Available: https://www.gurobi.com